\documentclass[utf8]{article} 

\usepackage[noblocks]{authblk}

\usepackage{url,lineno,microtype,subcaption}
\usepackage[onehalfspacing]{setspace}

\usepackage{amsmath,amssymb}
\usepackage{graphicx,color} 
\usepackage{array}
\usepackage{multirow}

\newcommand{\bs}{\boldsymbol}
\newcommand{\E}{\mathrm{E}}
\newcommand{\Var}{\mathrm{Var}}

\newcommand{\dd}{\mathrm{d}}

\author[1]{Marco Favino}
\author[1]{Alessio Quaglino} 
\author[1]{Sonia Pozzi} 
\author[1]{Igor V. Pivkin}
\author[1]{Rolf Krause}
\affil[1]{Institute of Computational Science, Universit\`a della Svizzera italiana}
\date{}

\title{Multiscale modeling, discretization, and\\ algorithms: a survey in biomechanics}

\providecommand{\keywords}[1]{\textbf{\textit{Keywords: }} #1}

\begin{document}

%

\maketitle

\begin{abstract}

Multiscale models allow for the treatment of complex phenomena involving different scales, such as remodeling and growth of tissues, muscular activation, and cardiac electrophysiology.
Numerous numerical approaches have been developed to simulate multiscale problems. However, compared to the well-established methods for classical problems, many questions have yet to be answered.
Here, we give an overview of existing models and methods, with particular emphasis on mechanical and bio-mechanical applications.
Moreover, we discuss state-of-the-art techniques for multilevel and multifidelity uncertainty quantification.
In particular, we focus on the similarities that can be found across multiscale models, discretizations, solvers, and statistical methods for uncertainty quantification.
Similarly to the current trend of removing the segregation between discretizations and solution methods in scientific computing,
we anticipate that the future of multiscale simulation will provide a closer interaction with also the models and the statistical methods.
This will yield better strategies for transferring the information across different scales and for a more seamless transition in selecting and adapting the level of details in the models.
Finally, we note that machine learning and Bayesian techniques have shown a promising capability to capture complex model dependencies and enrich the results with statistical information;
therefore, they can complement traditional physics-based and numerical analysis approaches.

\end{abstract}


\section{Introduction}\label{sec:intro}

Computational Science (CS), i.e., the study and the development of mathematical models, simulation methods, and solution algorithms, has become indispensable
for research and experimentation in science, engineering, medicine, and technology, cf.~\cite{PITAC_2005}.
CS allows to replace physical and laboratory experiments by simulations, i.e., \lq\lq virtual experiments\rq\rq, in order to predict the response of a system or optimize its design.

A mathematical model is usually derived from first principles and gives rise to initial and boundary value problems involving systems of, respectively, Ordinary Differential Equations (ODEs) or Partial Differential Equations (PDEs).
From a single-scale and single-physics point of view, current simulation methods are quite elaborated.
In fact, the existing methods for solving ODEs and PDEs have reached a high level of maturity: a-posteriori error estimators, adaptivity, sensitivity analysis, and fast solvers can guarantee accuracy and efficiency for a large class of problems.
Classical simulation approaches like Finite Element (FE) methods are routinely used in research and industry, and computational methods like Molecular Dynamics (MD) allow for a detailed investigation of effects at small scales, for example for the design of new materials and drugs.
However, many of the underlying mathematical models have been simplified to single-scale problems, with the idea of reducing the complexity of the arising equations
and make them affordable with the computational power of computers.

With the ever growing computational capabilities, the restrictions in terms of computational time have become less important.
As of 2017, available supercomputers deliver up to $93$ petaFLOPs \cite{TOP500} and exaFLOP machines are expected in the near future.
Thus, the available computational power allows to face even demanding simulations of complex and coupled models involving different length and time scales, i.e., {\em multiscale models}.

In parallel to the growth of computational capabilities, within the last decades numerical techniques for single-physics problems have been adapted for simulating complex systems and new multiscale numerical approaches have been developed.
In particular, we can name micro-Finite Element ($\mu$FE) method, variational multiscale (VMS) methods, and particle methods  in the context of \emph{multiscale discretization methods}, and
multigrid, multilevel, and cascadic methods in the context of \emph{multiscale solution methods}.
These techniques have been employed in different biomedical applications, such as remodeling and growth of tissues
(e.g., wound healing, bone healing, tumor growth), modeling of cardiovascular system, muscular activation, and electrophysiology.

The use of multiscale models and coupled systems characterized by a large number of parameters has also been a driving factor for the development of the field of uncertainty quantification (UQ),
aiming at evaluating the influence of the variation of such parameters.
This is of particular relevance for biological and bio-mechanical applications, where relevant parameters and the computational domain are not exactly known.
In this new perspective, we are moving from a deterministic approach of computing a single solution or a single trajectory, assumed to be uniquely determined by the initial data and the chosen parameters, to a probabilistic view.
However, this transition also requires an enormous amount of additional computational power. In this context, \emph{multiscale stochastic methods}, such as multilevel and multifidelity Monte Carlo,
are yet again playing a fundamental role for making these simulations affordable.

\begin{table}[hbt]
\begin{center}
\begin{tabular}{|c|c|c|c|c|}
\hline
                                & Model & Discretization & Solver & UQ \\
                                & (Section 2) & (Section 3) & (Section 4) & (Section 5) \\
                                \hline
$e<\epsilon$            & $e = \| u-\bar{u} \| $ & $e=\| u - u_h \|$ & $e=\| u^k_h - u_h \|$ & $e=\E[\|\hat{u}_h - \E[u]\|^2]$ \\
                                &&& & \\
                                & Simplified physics & Adaptivity & Exact coarse & All model \\
  $\Downarrow$          & Surrogate & Space enrichment & Approx. coarse & All discretization  \\
                                & Projection & Averaging & Projection & All solver  \\
                                &&& & \\
 Operation              & PDE creation & PDE solution & Smoothing & Quadrature \\
                                &&& & \\
                                & Sequential & Space decomposition  & Multilevel sum & $\alpha=1$ (MLMC) \\
 $\Uparrow$             & Concurrent & & Line search & $\alpha \ne 1$ (MFMC) \\
                                & Embedded && & Bayesian regression \\
\hline
\end{tabular}
\caption{Summary of all of the multiscale approaches treated in this work. We denoted by $\Downarrow$ the generation of a coarser scale, by $\Uparrow$ the coupling across different scales, and by $\bar{u}$ the true solution of the full-scale physical system.}
\label{table:summary}
\end{center}
\end{table}

In table \ref{table:summary}, we summarize the salient aspects of the multiscale approaches treated in this work. In particular, we identified four fundamental building blocks of multiscale methods:
\begin{enumerate}
        \item The quality measure $e$. This is also known as the discretization error in numerical analysis or the mean squared error in statistics. It represents the quantity that we wish to minimize for a given computational effort.
        \item The procedures from decreasing the resolution (denoted by $\Downarrow$). These are typically referred to as coarsening or model reduction. The former has enjoyed a widespread popularity, since it can change the resolution without the burden of creating of new model. However, it has a limited potential to reduce the computational effort and it cannot be applied in general. The latter is a truly black box technique and is nowadays enjoying a surge in popularity thanks to the rise of ever more capable machine learning algorithms. However, as such, it is limited to low- to moderate-dimensional systems and it produces results which are not easily interpretable from the physical standpoint. More details are discussed in Section 5.
        \item The core operation performed at each scale, which is also referred to as level in the context of solvers and UQ. Each level represents a single-scale and single-physics problem and as such can be dealt with using standard techniques.
        \item The coupling across different scales (denoted by $\Uparrow$). This is also known as interpolation in the context of multigrid methods or as control variate or variance reduction in that of UQ. It can be noted here that there are similarities between the approaches, e.g., between the multilevel sum of multigrid methods (Section 4) and of multilevel Monte Carlo (Section 5), as well as multigrid with line search and multifidelity Monte Carlo. However, there are also specific approaches, such as Bayesian regressions, that do not (yet) have similar counterparts in more than one area.
\end{enumerate}

In the context of high-performance computing software, a recent trend has been to recognize the limitations posed by the traditional segregation of numerical libraries, such as those implementing FE assembly routines and linear solvers. It has been shown that by overcoming this limitation and intertwining the two classes of methods, more efficient scientific software can be produced \cite{badia2017fempar}. In a similar way, we argue that the segregation of multiscale techniques for models, discretizations, solvers, and UQ, can be surpassed. This is true in particular in the context of the $\Downarrow$ and $\Uparrow$ operations discussed above and in table \ref{table:summary}.

With the advent of more advanced statistical techniques, such as deep learning, it is now possible to generate better and faster reduced models. Moreover, it is possible to shift paradigm from the traditional one of scale {\em selection} to that of scale {\em fusion}. In this new paradigm, the availability of multiple scales (e.g. models) is not seen as a hurdle, but rather as a strength, since each scale might contribute to explain part of the information filtered by the other ones. The pivotal role in this change of perspective is  currently being played by UQ.

In this paper, we present an overview on existing multiscale models and methods, with particular emphasis on mechanical and bio-mechanical applications. Moreover, we discuss state-of-the-art techniques for scale fusion in UQ. Finally, we argue that a higher level of abstraction could benefit all fields of multiscale methods, in particular for the improvement of the $\Downarrow$ and $\Uparrow$ operations.

\section{Multiscale Modeling}
 \label{sec:param-transfer}

Multiscale modeling refers to different approaches to derive equations which couple micro- and macro-scales.
The development of coupled multiscale models has been starting about 25-30 years ago in different fields.
We refer to the survey~\cite{DGroen_etAl_2014}, where an overview by scientific communities can be found.
There, Astrophysics, Biology, Energy, Engineering, Environmental and Material Sciences are identified as major fields for multiscale simulations.
Additional surveys on multiscale models can be found in~\cite{CM_03,PL_04,BCC_04,ELV_04}.
For what concern survey in biomedical biology, we can name
\cite{southern2008multi}, \cite{walpole2013multiscale}, and \cite{sloot2009multi} for multiscale modeling in biology.
For the particular case of biomechanics, we can name
\cite{deisboeck2010multiscale} and \cite{cristini2010multiscale} for multiscale modeling of tumors.
\cite{viceconti2012multiscale} for the mutiscale modeling of skeletal system,

In this section we provide a classification of the different techniques and keywords that are usually related to multiscale modeling with some examples.
A first general classification comes from the different approaches used to couple different scales:
\begin{enumerate}
\item in the \emph{hierarchical} or \emph{sequential} approach information is passed from the micro-scale to the macro-scale by means of \emph{homogenisation} techniques ~\cite{BLL_02,AG_05};
\item in the \emph{concurrent} approach different mathematical models coexist and are coupled, hence
requiring a continuous transfer of information between the scales.
\end{enumerate}
The second approach formally refers to applications in which micro- and macro-scale models are employed on different domains, such as in the coupling MD and continuum mechanics,
or they coexist on the same domain.
In the present classification, we prefer to employ concurrent for the former case,
and employ the term
\begin{enumerate}
\item[3.] \emph{embedded} approach for the cases in which the micro-scale is incorporated in the macro-scale model acting as a \lq\lq driving force\rq\rq~for this latter.
\end{enumerate}

\subsection{Sequential Multiscale Approach}

In the sequential multiscale approach, the scale transfer is realised by means of parameters that can be obtained from different strategies.
In contrast to the concurrent multiscale models, here the information transfer is unidirectional, i.e. only from micro- to macro-scale.
Homogenisation usually refers to methods that start from microscopical considerations
and then, by means of averaging procedures, equations or problem
parameters are upscaled to a coarser level~\cite{gloria:hal-00766743,NCharalambakis_2010}.
This process is realized by means of different techniques:
1) \emph{asymptotic homogenisation} that is a mathematical technique to homogenise highly oscillating parameters;
2) \emph{upscaling} in which detailed equations are solved at the micro-scale to obtain macroscopic information;
3) \emph{mixture theory} in which macroscopic equations are solved for all the constituents on the same domain
and coupling terms are employed to model the interaction at the microscale.

\subsubsection{Asymptotic homogenisation}

This technique is employed to study PDEs with highly oscillating coefficients which represent different material properties in the domain.
Diffusion-like equations are usually written as

$$-\nabla \cdot a^{\epsilon} \nabla u^\epsilon = f,$$
where $a^\epsilon$ denotes an oscillating parameter at the spatial scale $\epsilon$.
The discretization of such an equation would require a mesh of typical size $\epsilon$ in order to catch all characteristics of the solution.
Asymptotic homogenisation techniques try to derive a macroscopical equation

$$-\nabla \cdot a \nabla u = f,$$
where the parameter $a$ represents the effect of $a^\epsilon$ at the macro-scale.

\subsubsection{Upscaling}

\begin{figure}[h!]
\begin{center}
\includegraphics[scale=0.3]{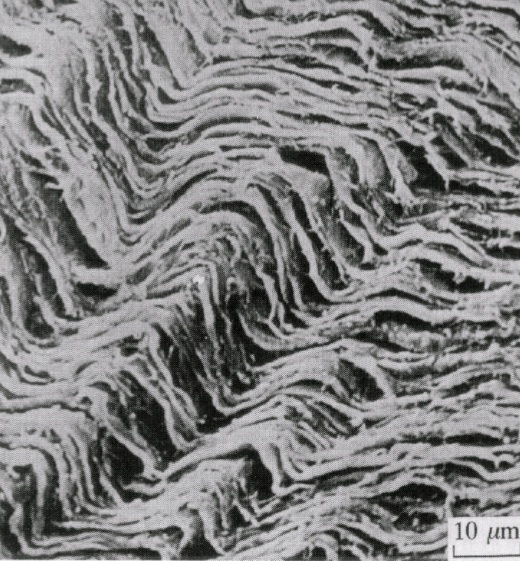}
\end{center}
\caption{Example of RVE for a biological tissue: a microscopic section of the periodontal ligament which is composed of two different phases.}
\label{fig:RVE}
\end{figure}

Upscaling techniques assume the existence of a Representative Volume Element (RVE),
which supplies all necessary information about the effect of the microscopical composition on the macro-scale.
RVE should contain a detailed structure of all the constituents of a heterogenous continuum.
We present the case of a two-phase material whose micro-structure is reported in Figure~\ref{fig:RVE}.

In the RVE, two different subsets $\Omega_s$ and $\Omega_f$ for solid and fluid phases are identified and different continuity equations
are formally solved on each of them

\begin{equation}
\begin{array}{c c c c}
\frac{\partial }{\partial t} \rho_s & = & \nabla \cdot J_s & \quad \text{on } \Omega_s\\[2mm]
\frac{\partial }{\partial t} \rho_f & = & \nabla \cdot J_f & \quad \text{on } \Omega_f\\[2mm]
J_s  \cdot n_s &  = & - J_f \cdot n_f & \quad \Gamma 
\end{array}
\label{eq:continuity}
\end{equation}
under the hypothesis that a separation surface $\Gamma$ between the constituents is known.
In Equation \eqref{eq:continuity} the unknown $\rho$ denotes a physical quantity object of the balance law, e.g. density, linear and angular momenta.
The microscopic solutions are hence upscaled by means of volume- or mass-averages~\cite{Hassanizadeh1989}
to obtain a system of equations at the macroscopical level.
In this process, also stochastic influences can be incorporated~\cite{ABourgeat_2004}.

Since this technique depends on the shape of the interface between the two materials,
effective analytical computations are feasible only for simple geometries.
For complicated geometries, an alternative technique is the numerical homogenisation.
This technique is widely used in rock mechanics, see e.g. \cite{hunziker2018seismic,hunzikerseismic}.
Here, the relevant macroscopical parameters are computed employing numerical simulations at the micro-scale on a RVE.
In~\cite{Cottereau201335} an iterative approach is described,
in which the material parameters obtained for a stochastic micro-structure are updated in an iterative manner within a coupled simulation.
Stochastic effects on the micro-scale are also considered in~\cite{XYin_etAl_2008},
where a statistical multiscale approach is presented, which aims at quantifying the influence of random material microstructure on material constitutive properties.

\subsubsection{Mixture theory}

Mixture theory assumes a material consisting of as many overlying continua as are the constituents in the system.
In this macroscopical approach, models are derived by imposing modified balance laws for all the constituents on the same domain $\Omega$.
They read

$$\frac{\partial }{\partial t} \rho_\alpha  =  \nabla \cdot J_\alpha + I_\alpha  \quad \text{on } \Omega,$$
where the subscript $\alpha$ refers to $\alpha$-th constituent.
The source terms $I_\alpha$ represent the interaction at the micro-scale between the different constituents,
and in case of balance they have to sum to zero.
In order to obtain a closed mathematical system, further hypotheses have to be introduced on the different components of the system.

\subsubsection{Poroelasticity Models for Biological Tissues}

Poroelasticity equations are an interesting example of models that can be formally derived from upscaling and mixture theory.
Very often the two approaches are believed to coincide but 
their equivalence is actually valid only in the linear regime~\cite{Schanz2003},
In fact, both derivations have their peculiarities.
Upscaled models allow for a detailed representation of the micro-scale \cite{Fritsch2006151} but their validity is limited to the linear regime.
On the other hand, models arising from mixture theory allow for a generalisation to the non-linear regime~\cite{biphasic2012},
but the micro-structure is taken into account only by means of constitutive laws.

In biomechanics, a standard single-phase description is not satisfactory to reproduce the mechanical response of biological tissues
as their behaviour strictly depends on their microstructure and composition.
In the recent years, poroelastic and multiphasic models are often more employed to obtain a more detailed description, e.g. for
cartilage \cite{Federico20052008},
intervertebral discs \cite{MalandrinoNL11},
heart walls \cite{Huyghe1991527}, and
periodontal ligament \cite{favino2011,favino2016nonlinear,FBK2014}.

\subsection{Concurrent Multiscale Approach} 
\label{sec:concurrent}
 
\begin{figure}[h]
\centerline{\includegraphics[width=0.3\textwidth]{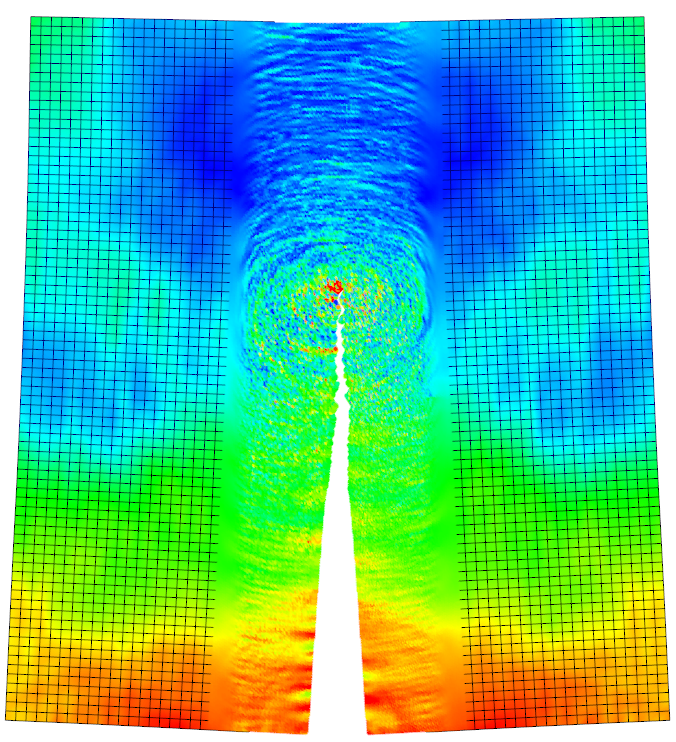}}
\caption{Coupling between Molecular Dynamics and Finite Elements (Mode-I fracture)\label{fig:mode-1-crack}.}
\end{figure}

In computational mechanics, the most prominent multiscale methods are designed for the coupling of MD on the micro-scale and continuum models on the macro-scale.
We refer to Figure \ref{fig:mode-1-crack} for an example from fracture mechanics illustrating this concept.
Here, a MD solution around the crack tip is combined with a FE simulation \lq\lq far away\rq\rq from the crack.
In these applications, the solution methods are based on a spatial decomposition of the
computational domain, which can be non-overlapping, partially
overlapping, or fully overlapping. Examples are the coupling of length
scales method, which uses a non-overlapping decomposition with a
lower-dimensional interface for coupling \cite{BABK_99}, the Bridging
Domain method  by Belytschko and Xiao, which uses a transition
domain for the information transfer between the scales, and the
bridging scale method by Liu et al.~\cite{WL_03}, which is designed
for the fully overlapping case but then reduces the size of the atomistic region by using absorbing boundary conditions outside the ``region of interest''. Recently, these methods have been extended to include also the randomness of the material structures using homogenization approaches, see~\cite{XYin_etAl_2008,Liu_etAl_2010} and the references cited therein.
In this context, also the Arlequin method has to be mentioned, which provides a general approach to multiscale coupling in particular in mechanics, see~\cite{HBDhia_1998}.
In addition to multiscale approaches where MD and finite elements are coupled, also other coupled models exist. An example can be found in~\cite{May_etAl_2011}, where finite elements are coupled with discrete automata for modeling tumors.

We finally note that sometimes also multi-dimensional models are referred to as multiscale models.
Here, $1$D (often ODEs) and $3$D models (often PDEs) are coupled,
see~\cite{SGrein_et_al_2014,TPasserini_etAl_2009,RKrause_JSteiner_2010}.
The fully overlapping case is described in more detail in Section~\ref{sec:embedded}.

\subsubsection{Scale Transfer} A particular challenge in the framework of coupled multiscale simulations is
the transfer of discrete fields - in mechanics usually displacements and stresses - between the different scales.
By design, the models on the different scales will have different characteristics.
It is well known that a simple and straightforward coupling of MD and finite elements
leads to spurious oscillations and/or to the so called ``ghost forces'' at the coupling interface, which can spoil the accuracy of the simulation.
Moreover, in MD the displacement of the atoms are  pointwise given quantities,
whereas in a continuum based model such as finite elements the computed displacements live in a suitable function space  (usually the Sobolev space $H^1$).
From neither a physical nor a  mathematical point of view it is  possible to identify finite element displacements with MD displacements directly.
Instead, a suitable transfer operator has to be derived which allows for a stable and seamless transfer of discrete fields between the different models.
See~\cite{FackeldeyKrauseKrauseLenzen09} and the references therein for details on this aspect.

A similar -but less severe- problem shows up with respect to the discretization in time.
Here, the different time scales of MD and continuum based approaches have to be taken into account.
One possibility is to use time integrators as  the SHAKE-RATTLE integrator~\cite{HLW_02}, which is a standard symplectic integrator for constrained Hamiltonian systems.

\subsection{Embedded Multiscale Approach}
\label{sec:embedded}

Embedded multiscale models are characterised by two sets of equations which describe the micro- and macro-scale defined on the same computational domain.
These equations are usually strongly coupled: the microscale is incorporated as a \lq\lq driving force\rq\rq~in the macro-scale equations.

Cardiac electrophysiology and cardiac electromechanics are classical examples of embedded approaches.
In the cardiac tissue, cells maintain a difference of ion concentration between the interior and exterior of the cell by means of active ionic (K$^+$ and Na$^+$) pumps.
Hence, the cell membrane acts as an insulator, keeping a difference between the intra- and extra-cellular potential (known as membrane voltage).
Since a cellular approach for the simulation of the heart would be unaffordable,
cardiac tissue is usually modelled by means of mathematical homogenisation \cite{Keener:1998:MP:289685,franzone2014mathematical}.
A macro-scale continuum approach is employed, describing the spatial concentration of the intra- and extra-cellular ionic charges.
This bi-phasic representation leads to a system of equations known as bidomain model whose components are 
\begin{itemize}
\item a set of Hodgkin and Huxley-like ODEs~\cite{Hodgkin28081952}
describing the dynamics of the ionic pumps and ionic concentrations at the micro-scale;
\item two reaction-diffusion PDEs modeling the propagation of intra- and extra-cellular potentials at the macro-scale.
\end{itemize}
The evolution of ionic pumps and concentrations is used to compute the macroscopic ionic current (i.e. the reaction term of the PDEs), providing the coupling from the cell scale to the macro-scale.
On the other hand, coefficients of Hodgkin-Huxley depend on the membrane voltage, providing the coupling from the macro-scale to the cell scale.
This bi-directional coupling between the different unknowns at the different scales is characteristic of embedded models.

In order to include stochastic effects in the parameters of ionic pumps, methods based on both Markov Chain Montecarlo and stochastic differential equations have been proposed~\cite{CDangerfield_etAl_2012.png}.
The former is in general computationally expensive and does not allow to include the simultaneous variation of different parameters.
The latter provides a more simple and efficient way of including uncertainties in the material parameters of Hodgkin-Huxley models.
The stochastic variability of parameters in the ionic pump ODEs have been shown to be of particular effect for isolated cells
but not to be relevant for a complete cardiac tissue \cite{EPueyo_etAl_2011}.

Growth and remodeling in continuum mechanics are other relevant examples of embedded approaches.
In this case, an inelastic deformation tensor is introduced to describe the evolution of the stress-free configuration \cite{biphasic2012}.
The evolution of this term is hence coupled with ODEs,
which describe the evolution and exchange of material with nutrients at the micro-scale.
Coefficients may depend on stress state of the material, again leading to a 
A last case of embedded approach is the statistical orientation of fibers inside cartilage~\cite{Federico20052008}.
Continuum mechanics equations are coupled with an evolution equation that governs the mean orientation of the fibres.

\section{Multiscale Discretization Approaches}

\begin{figure}[h]
\centerline{\includegraphics[width=0.4\textwidth]{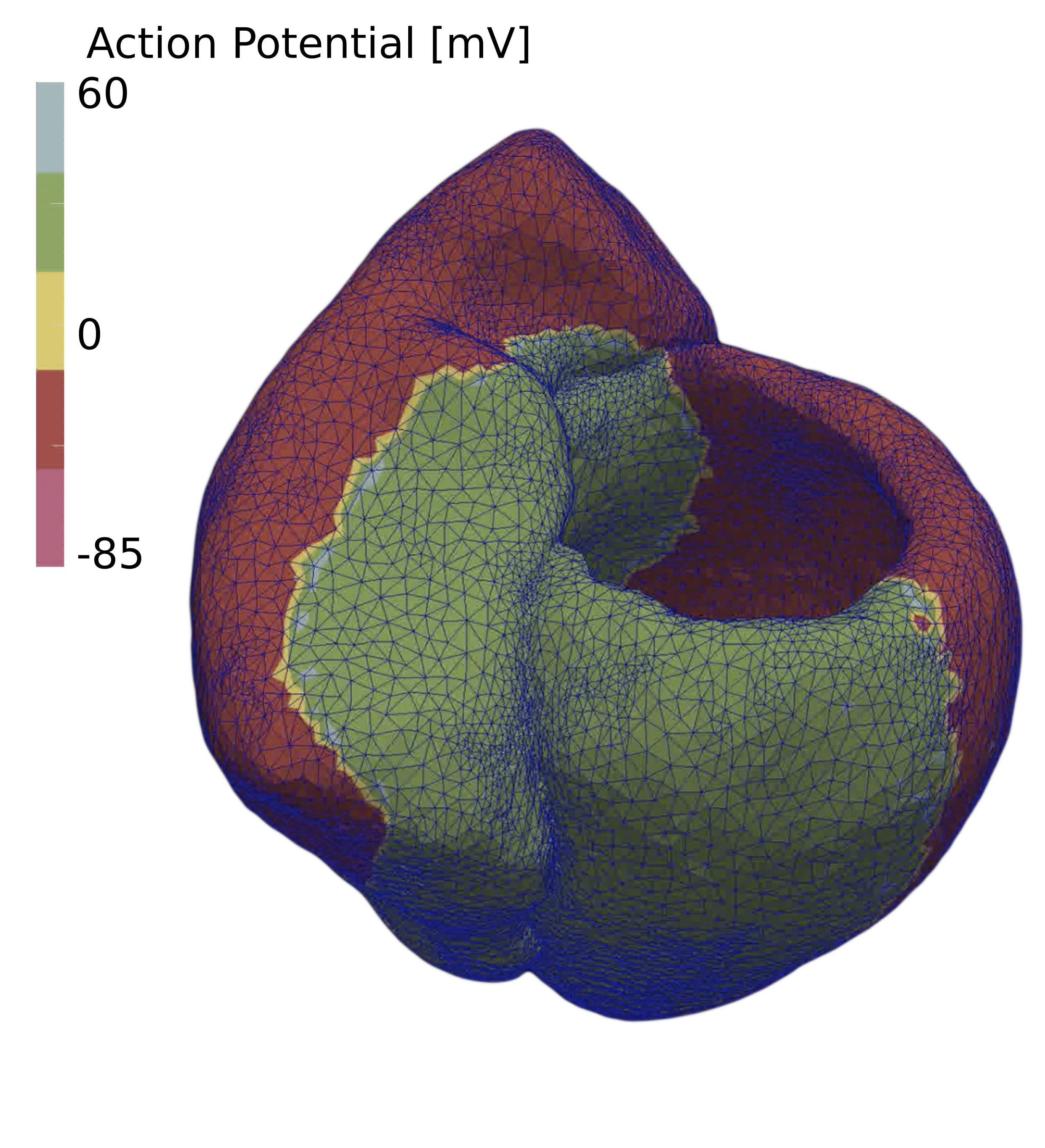}}
\caption{Propagation of the potential in a cardiac tissue: steep space-time gradients exist requiring a fine mesh and small time-steps in order to correctly them.}
\label{fig:activation}
\end{figure}

Solutions of equations arising from multiscale models still manifest their multiscale origin being characterized from steep gradients and large time derivatives.
For example this can be observed in the bi- and monodomain model for cardiac electrophysiology.
An example of solution is reported in Figure~\ref{fig:activation}.
Here we can observe large zones at the resting potential and large activated zones which are connected by a small region where the wave front is localized.
This behaviour is typical of reaction-diffusion equations.

Also consolidation problems of poroelastic biomaterials are characterized by a very rapid increase of pressure and then a relaxation of the tissues.
In order to reproduce this very quick transient, a very short time step is needed
and hence a large sequence of non-linear implicit problems has to be solved.
The short time-step also imposes a very strict condition on the mesh size \cite{doi:10.1061/9780784412992.110}.

In order to correctly reproduce catch these boundary layers in the discretization of such models,
multiscale discretization methods, which are object of this section, have been developed.

High resolution approaches, such as $\mu$-FE, have been developed to catch small scale effects but in general they results in large systems of equation.
Hence, other techniques, such as VMS, have also been developed allowing to incorporate small scale effects into coarser discretization combining
analytical information on the solution with standard FE methods.
Ideas similar to this latter are also used in exponential time integrator for stiff ODEs.
Here the linear part of the equation for which an analytical solution is available is exactly integrated while numerical methods are employed for the high order terms.

A totally different approach comes instead from particle discretizations usually employed in fluid-dynamics and MD.
These methods allow to remove the strict condition given by the CFL condition and do not require the creation of a mesh.
Coarse-grained discretizations allow again with a multiscale approach to reduce the computational burden given by particle approaches creating aggregates of particles.

\subsection{High Resolution Approaches} 
The possibly most straightforward approach to  treating different length-scales in bio-mechanics might be the usage of high resolution models based on single-scale models.
Here,  prominent example is the simulation of trabecular bone using the so called micro-finite element method:
a high resolution finite element mesh is used which allows for resolving the structure of the trabeculae.
Naturally, this approach finds it limits where the material does not behave as a continuum any more, as this is a necessary prerequisite for the application of continuum based material models and the finite element method. The resulting large scale systems are then solved using established and well scaling fast iterative methods as multigrid methods on massively parallel machines~\cite{PArbenz_etAl_2008}. In the case of non-linear material models, the situation is more complex,
but it can be handled efficiently even for the arising large scale systems, cf., e.g.,~\cite{DChristen_etAl_2013}. 
These high resolution models, which have been originally designed for simulating the mechanical behaviour of the bone, have also been extended in order to investigate, e.g., bone-cell response to mechanical stimuli~\cite{Lacroix20065326}. In~\cite{FGerhard_etAl_2009} it is argued that, in order to model and simulate this interaction more accurately,
hierchical models should be developed and employed.
An approach to stochastic modeling can be found in~\cite{Thomsen:1994gx}.
Here, stoachastic effects are added to a continuum based model.

\subsection{Coarse Graining}
A fundamental problem for numerical simulations on the molecular scale is the fact that the currently feasible time- and length scales are still too small for many practical applications.
Although standard simulations tools in MD show excellent scaling behaviour also on larger parallel machines,
the computation times needed for simulations in classical MD are still very high:
an example might be  the simulation of a complex virus in a time range of about a milli-second.
In order to address this problem, so called coarse grained (CG) methods have been introduced,
which aim at using coarser scales in space and time without sacrificing too much the microscopic information, see, e.g.,~\cite{Voth09,niels04}.
Although there is a huge variety of CG models available in the literature, they share the common idea of partitioning the system into single interaction centres, the so called \lq\lq beads\rq\rq.
For this in-space-reduced system, the CG potentials are derived, which provide the effective force field for the interaction of the CG particles or beads.
Often, these CG potentials are given in a parametrised form and they aim at preserving structural properties observed in atomistic simulations.
Examples are approximations based on the potential of mean force~\cite{shell01,mue00}
or the inverse Monte Carlo technique~\cite{murto04}
and Inverse Boltzmann.
CG potentials might also be designed to preserve thermodynamic properties~\cite{shell01,mar04,mar07,mar08}. 

Multiscale Coarse-Graining (MSCG)~\cite{izve05} is a different class of CG potentials that 
uses information from high-resolution simulations or models for building the CG potential.
This is usually done by means of a force matching method~\cite{par04,erc94}.
For the analysis of MSCG in the framework of statistical mechanics, see \cite{mscg01,mscg02}.
From a numerical point of view, the computational effort for CG simulations can be reduced by choosing a new multi-resolution basis~\cite{dasa12,dasb12}.
With respect to stochastic approaches, combinations of CG approaches with Monte-Carlo methods have been derived and analyzed, see~\cite{AKatsoulakis_et_Al_2003}.

In terms of biological and biomechanical applications,
CG models have been used for simulating the behaviour of platelet aggregation as a stochastic process~\cite{IPivkin_PRichardson_GKarniadakis_2014}.
Also on the cellular level, in~\cite{FVermolen_AGefen_2012} an approach is presented
which incorporates  cell death and proliferation as stochastic processes.
A multiscale description for wound healing using different models on cellular, cell colony, and tissue scale is described in~\cite{FVermolen_AGefen_2013}.

\subsection{Dissipative Particle Dynamics}

Dissipative particle dynamics (DPD) is a relatively new, potentially very effective particle-based approach, which was applied in simulations of physical and biological systems from atomistic to macroscopic size.
The DPD model consists of particles (\lq\lq beads\rq\rq) which correspond to coarse-grained entities, representing clusters of atoms or molecules~\cite{Hoogerbrugge_SMH_1992,Pivkin_DPD_2010}.
The size of these clusters, which defines the level of coarse-graining in DPD simulations, can vary by many orders of magnitude in different applications.
For simple fluids, the DPD particles interact with each other via pairwise forces, e.g., conservative, dissipative, and random~\cite{Groot_DPD_1997}.
All forces in DPD conserve linear and angular momenta.
Conservative force in DPD is derived from soft quadratic potential and its magnitude depends on the physical system.
Dissipative force takes into account dissipation of energy due to friction forces.
Stochastic component of the DPD, the random force, takes into account the degrees of freedom which were eliminated as a result of coarse-graining.
Dissipative and random forces form DPD thermostat~\cite{Espanol_SMO_1995}.
In simulations of complex systems other forces are usually added, which include bonded interactions, angle potentials, etc.
Unlike MD, where the choice of forces is based on a theoretical model of the physical system to be simulated, the DPD model involves forces of a form independent on the physical system.
Therefore, parameters in DPD simulations have to be carefully chosen in applications.

Two interpretations of the DPD method exist.
According to the first one, DPD can be considered as a stochastic CG molecular dynamics~\cite{Groot_DPD_1997}.
Here, each DPD particle corresponds to a cluster of a small number of atoms or molecules, typically in the range between 3 and 10.
In this case, the most common process of choosing DPD simulation parameters for fluid systems is based on assuming that the dimensionless compressibility of the DPD fluid is equal to the dimensionless compressibility of the physical fluid~\cite{Groot_DPD_1997,Keaveny_ACS_2005}.
The time scale is defined by matching the diffusion constant of the molecule of interest~\cite{Groot_MSO_2000,Groot_MSO_2001} or vorticity~\cite{Keaveny_ACS_2005}.
The level of coarse-graining in this interpretation of the DPD is limited however~\cite{Trofimov_TCI_2003,Pivkin_CGL_2006}.
For systems containing water it is usually less than 10 water molecules per single DPD particle.
Nevertheless, due to the soft potentials employed in DPD, the time scales are typically quite large comparing to MD, resulting in speed-up factors of more than $10^4$~\cite{Groot_DPD_1997}.

The second interpretation of the DPD method is based on the fact that DPD forces conserve linear and angular momenta.
As a result, DPD provides correct description of the hydrodynamic interactions in simulations even with relatively small number of particles~\cite{Espanol_HFD_1995}.
Here, each DPD particle can correspond to a cluster of a very large number of atoms or molecules.
The DPD simulation parameters are chosen so that relevant dimensionless numbers (such as Reynolds number in flow simulations) or specific properties of fluid are preserved.
Thermal fluctuations are still present in simulations, and therefore the method is typically used for modeling mesoscale systems.

The dual nature of the DPD method, which allows simulations of atomistic as well as mesoscopic systems, can be effectively utilized in multiscale modeling approaches.
In coupling of MD and continuum based methods (e.g., finite element), DPD can be used as an intermediate model, helping to bridge the gap between scales.
On the other hand, coupling DPD itself with continuum solver will allow simulations of systems, in which the smallest resolved length and time scales will be controlled by the level of coarse-graining used in DPD.
Uncertainties and unresolved small scale effects in the system then can be integrated into the stochastic nature of the DPD method.
The difficulty and the effort required for such coupling should be similar to the coupling of MD and finite element method.
However, due to the larger time scales in DPD, the hybrid multiscale system based on DPD will allow for larger simulation times as compared to MD based hybrid systems.

DPD has been applied to model many complex multiscale systems, including simulations of polymer solutions, brushes and melts~\cite{Schlijper_CSO_1995,Malfreyt_DPD_2000,Spenley_SLF_2000}, binary mixtures~\cite{Coveney_CSO_1996}, amphiphilic systems~\cite{Jury_SOA_1999,Shillcock_ESA_2002,Peter_APC_2014}, and cells~\cite{IPivkin_GKarniadakis_2008,Peng_LBA_2013}.
Also, certain processes on the cellular level were modeled stochastically using multiscale approaches in order to allow for a realistic yet simple description of complex behavior, for example see~\cite{Pivkin_EOR_2009}.

\subsection{Exponential time integrators for stiff ODEs}

The simulation of the monodomain equation is usually realised by means of semi-implicit methods~\cite{Franzone,Ethier,Veneziani}.
Often, the non-linear term is treated in an explicit way and the  diffusive term implicitly.
This allows to employ explicit time integration for the ODEs describing ionic concentrations, thus realizing the multiscale coupling. High order semi-implicit schemes can be formally derived, but the strict stability condition depending on the mesh-size prevents their use in realistic simulations. To avoid the usage of very small time-steps, Backward Differentiation Formula can be used~\cite{Ethier,Veneziani,hundsdorfer2013numerical}. Unfortunately they require the expensive evaluation of the derivatives of the ionic currents.
Also, operator-splitting methods are commonly used to avoid the implicit evaluation of derivatives~\cite{Qu,Sundnes2}.
 In~\cite{campos2013comparing} different techniques (adaptive time step-methods, partial evaluation and lookup tables, and the exploitation of the code concurrency via OpenMP directives) to automatically speed up the numerical solution of cardiac models are proposed.

The solution of monodomain equation is characterised by a travelling peak with steep gradient around the so called depolarisation region.
It is localised in a relatively small zone while in the rest of the domain the solution is smooth.
Capturing this travelling wave renders the simulation of the electrophysiology of the heart  challenging due to a required high spatial resolution.
In~\cite{KrausePotseDickopf}, an efficient adaptive strategy for large parallel machines has been developed. However, the refinement is controlled using error indicators, as reliable and efficient error estimators for this kind of problem are not yet available.

A first step towards the development of error estimators for this kind of coupled multiscale problem may be the single-scale continuous formulation of electrophysiology recently presented in~\cite{DHurtado_DHenao_2014}.
There, a variational formulation of cardiac electrophysiology is given which allows for a completely continuous formulation of the underlying multiscale problem in terms of a minmax-problem.  
This formulation allows for deriving bounds for the time-steps below which the objective function is strictly convex, as it is highly beneficial from a numerical point of view.

\begin{figure}[h]
\centerline{\includegraphics[width=0.7\textwidth]{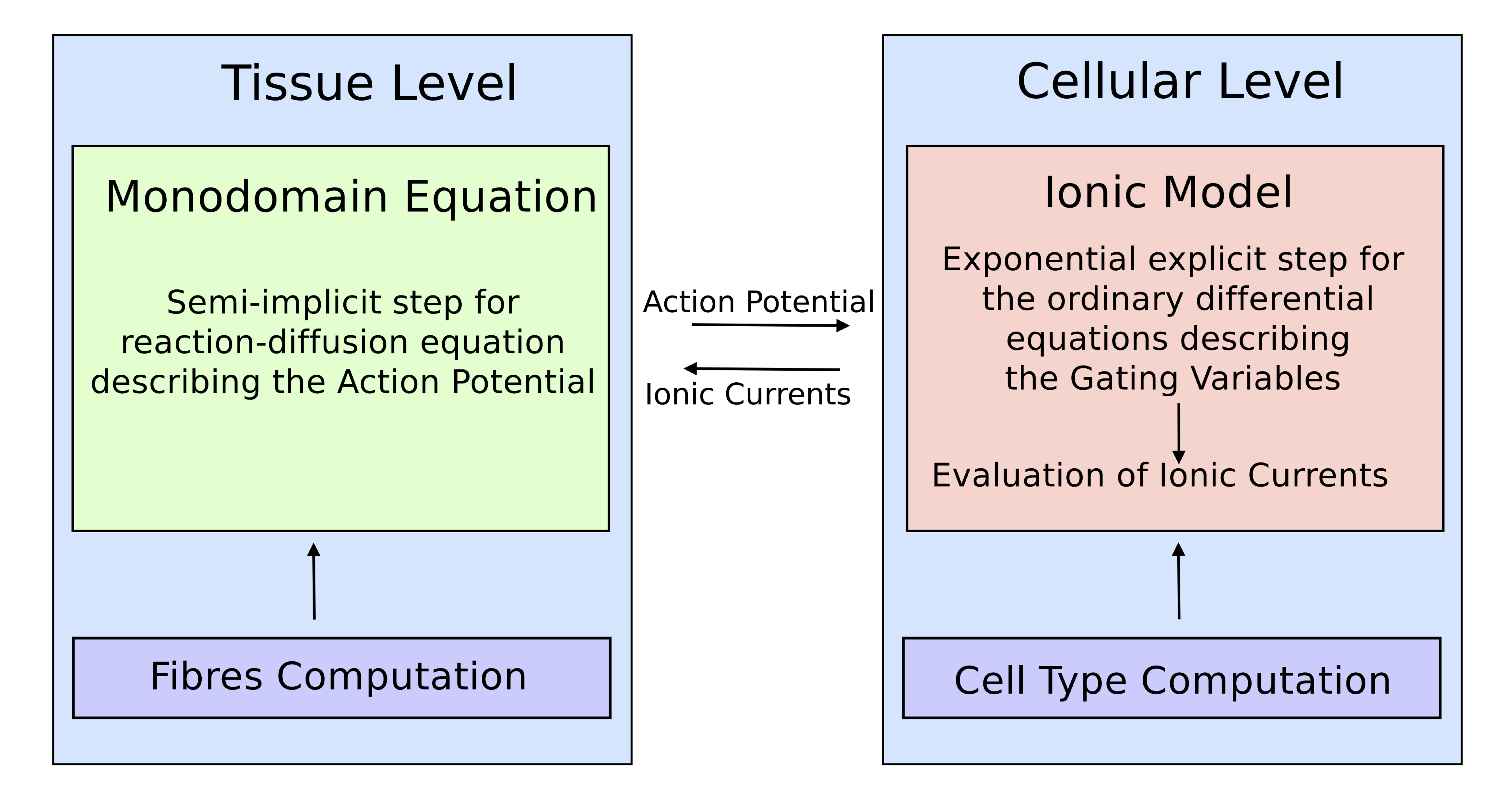}}
\caption{Schematic representation of the solution the multiscale solution strategy.}
\label{fig:data-exchange}
\end{figure}

At the cellular level, the main numerical difficulty in cardiac electrophysiology is the solution of the time-evolution of the gating variables. The Rush-Larsen (RL) scheme~\cite{RL} allows to overcome the above-mentioned problematic. In this approach, the time integration of the stiff ODEs describing the gating variables is performed assuming the action potential constant over the interval of integration.
This method is equivalent to an explicit exponential integrator~\cite{Norsett}
and prevents over- and under-shooting in the numerical solution of gating variables, ensuring that these belong to a physiological range. 
RL method improves the stability with respect to standard explicit integrators, but its convergence is limited to the first order.
RL scheme falls into the class of exponential integrators and it originates from a peculiar linearization of the original problem.

Literature of exponential integrators is well established~\cite{Lawson,Beylkin, Boyd, Cox, Smith,Hochbruck,Minchev,Krogstad}
but it is focused mainly on semilinear ODEs with constant coefficients of the linear terms.
When considering the ODEs describing the gating variable dynamic, the coefficients of the linear terms are not constant, rather  they depend on the action potential.

A recent attempt to derive a second order RL scheme has been proposed in~\cite{Sundnes}.
This method consists in a predictor-corrector middle-point method.
The predictor step exploits the standard RL method,
then, in the fashion of Heun integrators, the standard first order formulae are evaluated at the mid-point of the time step.
This approach has been shown to outperform standard RL and also Runge-Kutta methods with a double computational cost for each time-step~\cite{Gomes2015}.
In~\cite{Veneziani}, formulae for generalized RL schemes have been presented
but they were shown not to be in general A-stable.
An alternative strategy to employ high order exponential integrators consists in transforming the original problem in an equivalent one in which the leading coefficient of the right-hand-side is constant on the discretization interval.
These high order methods have been derived for ionic models in~\cite{Lontsi}, where a comparison with high order RL is performed.

\section{Multiscale Solvers}

In multiscale applications, small spatial and temporal effects interact larger scales.
Hence, small scales impose severe restriction on time discretization schemes and on mesh resolution in order to obtain physically meaningful results.
For this reason, in a standard Rothe's or line simulation framework,
these lead to several large linear systems to be solved at each time step of a time advancing scheme.
For example, in cardiac electrophysiology small time step-sizes are necessary for an accurate description of ionic channels,
while fine meshes are necessary at least where the wave front is localized.
In order to do this, adaptive mesh refinement approaches allow for a local refinement where numerical error is larger
but they are time-consuming in the process of error estimation and remeshing.
The use of uniform fine meshes on the other hand removes this computational burden
but requires efficient numerical methods for its solution.

Again, the idea to accelerate the solution of linear systems comes from multiscale decomposition.
Multigrid (MG) methods allows to obtain optimal complexity in the solution of the arising linear systems thanks to the use of different mesh refinements.
This is realized since the different meshes are able to reproduce different frequencies of the solution.

The first MG scheme developed for the solution of the Poisson problem was developed in the sixties by \cite{Fedorenko}.
The strong efficiency of MG strategies became more clear ten years later when \cite{Brandt} and \cite{Hackbusch} extended the MG idea to non-linear problems.

Convergence theory of MG methods is well-developed for symmetric and elliptic operators \cite{Braess,Briggs} employing the natural norm of the problem.
Nevertheless they showed good convergence rates and optimal behavior also for more general problem such as diffusion-convection,
saddle-point, and constrained problems.
In particular for convex quadratic functionals subject to box constraints such as the case of contact problems in computational mechanics
a monotone decrease of energy has been shown~\cite{Krause}.

In the implementation of MG strategies, two are the main instruments: restriction and projection operators, which allow to transfer information between one mesh and the other,
and smoothing operators, which allow to smooth the error in the solution when a correction is transferred from a coarse level to a finer one.
Nowadays, standard MG methods are often referred to as Geometric Multigrid Methods.
There are in general characterized by nested meshes,
standard FE assembly of the operators on each levels,
and restriction operators which coincide with the transpose of the interpolation operators.
Although this approach is the most efficient one both in the construction of the operators and in convergence rate,
it presents the severe limitation that the finer meshes have to be obtained from uniform refinement of a coarse one.

To overcome this problematic, several extension have been presented.
Algebraic Multigrid (AMG) methods allows to avoid any geometrical information on the coarse level by only employing the structure of the stiffness matrix.
Using aggregation methods they construct coarse levels and restriction operators.
AMG implementations, such as hypre/BoomerAMG, are available and allows for large parallel simulations~\cite{Henson}.
The main limitation of AMG strategies is that the coarsening strategies depends on the structure of the stiffness matrix, i.e. it has to be an M-matrix.
Hence, the underlying problem has to be scalar and constructed from low order discretizations.
Moreover, also reaction- and convection-dominated problems does not fall in this class unless multiscale or stabilized discretizations, such as VMS, are employed.

A suitable construction of the restriction operators is crucial for the efficiency of MG.
Semi-geometric multigrid (sGMG) is a recent technique for the use non-nested coarse levels.
The main idea behind is that starting from a fine mesh,
restriction operators can be constructed as the transpose of the projection operators from a coarse to a fine mesh.
This strategy seems promising, providing good convergence rates with the price of assembling the projection operator.
For the construction of coarse level operators, both standard FE assembly and Galerkin assembly can be adopted.
Any parallel linear algebra library (e.g. PETSc or TRILINOS) nowadays provide the framework for sGMG.

Stable and efficient techniques for the transfer of information between discrete fields on non-matching volume or surface meshes is also an essential ingredient
for the simulation of coupled multi-physics problems. In particular when the different \lq\lq physics\rq\rq~requires different mesh resolutions.
 From an High Performance Computing point of view,  the freedom to handle the different levels of refinement in a completely arbitrary way makes it possible to easily provide better balanced computations.
Moreover, considering not nested grid hierarchies it is possible to represent the considered geometries with the required level of accuracy.

Appropriate smoothers have also to be chosen and adapted for theproblem at hand,
since the standard receipt of $3$ Gauss-Seidel steps is not always effective.
This holds particularly true for vector and saddle-point problems where the use of point-block and Vanka smoothers is necessary to obtain reasonable convergence rates.

Finally, we want to point out that the MG idea can be generalized in several ways.
Instead of using coarse spaces arising from FE discretizations on coarse meshes, any coarse level which enjoy the approximation property can be used.
Another generalization may also come from the introduction of line-search techniques when the correction from the coarse level is added to the current solution \cite{livne2012lean}.

\subsection{Multiscale  in High Performance Computing}
With respect to the realisation of multiscale approaches in terms of simulation software, the main difficulty is arguably the fact that most of the existing simulation softwares have been written with a single-physics or single-scale application in mind. Thus,  exploiting the available computational power in the context of multiscale simulations turns out to be  difficult, as specially tailored algorithms have to be designed and implemented, which are capable of dealing with the massive parallelism of current and upcoming supercomputers.  This holds in particularly true if  different simulation methods - like MD for the micro-scale and finite elements for the macro-scale - have to be intertwined. As a matter of fact, their efficient simultaneous usage within a common simulation framework is far from trivial, see, e.g.,~\cite{DBLP:conf/eScience/KrauseK11}.

In fact, many multiscale simulations are usually carried out by combining different codes through more or less elaborated interfaces.
We refer to ~\cite{DGroen_etAl_2013}. Following this approach, data has to be transferred between the different scales (or simulation codes), see Figure \ref{fig:data-exchange}.
\begin{figure}[h]
\centerline{\includegraphics[width=0.7\textwidth]{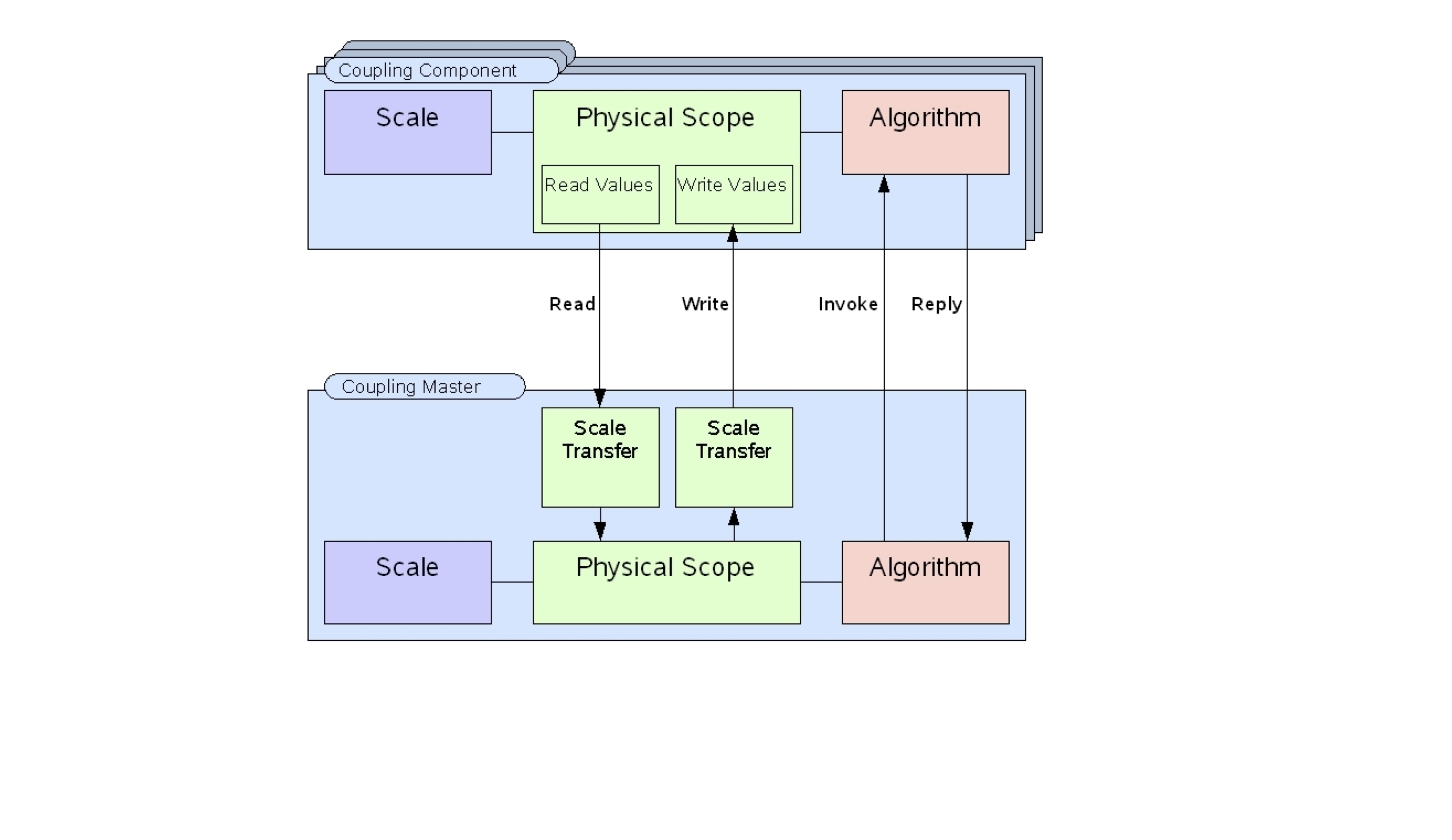}}
\caption{Data exchange between different simulation tools.\label{fig:data-exchange}}
\end{figure}
This is usually done in a sequential way, as the progress of the simulation on one scale will depend on the progress of the simulation on a different scale. One possibility to improve the computational speed is to advance the faster scale independently from the slower scale by doing a multi-rate time integration, cf., e.g.,~\cite{VSavcenco_RMMMattheij_2008}.
With respect to scalability on massively parallel machines, however, the inherently sequential approach used for the realisation of many multiscale simulations can limit dramatically the scalability. For example, for the coupling of MD and finite elements, in~\cite{DBLP:conf/eScience/KrauseK11} two parallel simulation tools have been coupled using a parallel scale transfer.
Nevertheless, the originally good scaling of the two simulation tools when used ``stand-alone'' could not be maintained. This is simply due to the fact that the respective computational demands of the simulations on the employed scales are substantiably different and, combined with the scale transfer, do not allow for a good load balancing.

Within the framework of MD and particle-based CG models, the situation is less severe, as both models are relatively similar. Thus, by now efficient softwares for MD simulations (GROMACS, LAMMPS) as well as approaches for their efficient coupling~\cite{doi:10.1021/ct400727q} exist. However, as explained above, this relatively good behavior can not be transferred to multiscale simulations including different models as, e.g., MD and finite elements.

When different models are intertwined nicely as for example in the case of cardiac electrophysiology, excellent scalability can be achieved on current supercomputer architectures, see, e.g.,~\cite{KrausePotseDickopf}.
This is possible since the smaller scale, i.e. the models for the ion-channels, can be  seamlessly integrated into the finite element method used for the diffusive part of the mono- or bi-domain equation.
The integration is done by simply solving a set of ODEs at each mesh node or at each quadrature point. As this additional work load is completely local and moreover scales nicely with the number of degrees of freedom for the finite element simulation, massively parallel machines can be used efficiently.

From the above examples it is clear that the possibility to use massively parallel supercomputers will heavily depend on the type of the multiscale model and on its characteristics. To put it simply, the better the different scales can be intertwined, the more likely good scalability will be. 

Making the computational power of large machines available for truly heterogeneous multiscale simulations thus remains a significant challenge. Coupling existing codes is possible, but can lead to a drastic loss in scalability. The ongoing trend in supercomputing towards heterogeneous hardware (CPU, accelerators, GPUs, shared memory, distributed memory,\ldots) and the corresponding diversification in terms of programming tools (OpenCL, CUDA, OpenMP, MPI,\ldots) will make it moreover necessary to  provide highly specific implementations. A possible remedy here could be to use asynchronous solution methods as stochastic solution methods, as those might allow for obtaining good scalability even for the strongly heterogeneous demands of multiscale simulations. \section{Multiscale Uncertainty Quantification}
\label{sec:UQ}

As discussed in the previous sections, simulations are characterized by uncertainty. UQ is a broad term, which can encompass a plethora of physical properties, as well as the lack of an infinite amount of data used for their estimation. However, it is useful to focus on some of its some properties, which can help to better characterize it:

\begin{itemize}
	\item Nature of the uncertainty. This is considered to be either aleatoric or epistemic \cite{smith2013uncertainty}. The former consists of all phenomena whose stochasticity is of interest and which, therefore, are naturally treated probabilistically. An example is the probability distribution of a certain parameter. The latter type stems instead from incomplete information, which can be due to limited data, noise, model errors, approximations, etc. While this uncertainty is not truly probabilistic, it is often conventient to model it as such, rather than use a better model (which is more expensive) or acquire more data (which contains noise).
	\item Sources of uncertainty. The uncertainty can arise from the parameters of a problem, from the model itself, or from its data. For example, material parameters, external loadings, boundary conditions, and geometry, are in practice not known exactly. Depending of the origin, different methodologies for treating the uncertainty might be applied. A survey which also discusses available simulation tools is presented in~\cite{GLin_etAl_2012}. 
	\item Direction of quantification. Essentially, most types of UQ can be classified as either forward uncertainty propagation or inverse uncertainty quantification. In the former case, the interest is to quantify the effect of the (input) uncertainties on some quantities-of-interest (QoI), which are typically referred to as sensitivities. This can be done independently from any data, by assuming a given probability distribution on the input. In the latter case, the aim is instead to solve and inference problem, starting from the data. The aim is to infer either the internal states of the model, to calibrate its parameters, or both.
	\item Dimensionality. In some applications, the uncertainties cannot be described as random variables but as random fields. This is the case for, e.g., spatially-dependent parameters. A possible method to approach this problem is to use the Karhunen-Lo\'eve (KL) theorem \cite{loeve1978}. The idea is to write a stochastic process  as the infinite sum of random variables - which can be interpreted as a stochastic Fourier series. Then, only a certain number of basis variables is employed to approximate the solution. Obviously, if the variability in the data is not large, few functions are sufficient to reproduce the random process. 
\end{itemize}

In this survey, we focus on forward propagation of high-dimensional uncertainties for models described by a stochastic PDE with random coefficients. We remark that stochastic PDEs become particularly challenging when the problem has also a multiscale structure. The multiscale nature arises when the material parameters have many small inclusions that are randomly distributed or when the covariance function has a sharp peak. In this case, the construction of a reduced basis requires a large number of terms. In order to characterize the influence of material properties on the micro-scale, various multiscale models for uncertainty quantification have been developed, cf.~\cite{MPapadrakis_GStefanou_2014,AChernatinskiy_etAl_2014}. Here, the stochastic models on the micro-scale are mainly used for characterizing unknown material properties. We finally note that a purely probabilistic model for material fatigue and formation of micro-cracks has been presented in~\cite{schmitz2013}.  The model is derived from reliability statistics and takes into account size effects, inhomogeneous strain, and thermal effect. Interestingly, this model acts as a "probabilistic post-processor'' on data created by deterministic macro-scale simulations.

\subsection{Quadrature methods}
The most straightforward approach in the context of uncertainty propagation is the computation
of integrals. If we denote by $\psi_h(\bs{s})$ the QoI, as computed using the
model for a value of the random input $\bs{s}$, then its mean
is obtained by:
\begin{equation}
\E[\Psi_h] = \int \psi_h(\bs{s}) ~p(\bs{s})\:\dd\bs{s}
\label{eq:expect}
\end{equation}
where $p(\bs{s})$ denotes the probability density of the random inputs.
Similarly, one can compute other statistical indicators of the QoIs
(e.g.moments, densities etc), by evaluating integrals of the form $\int
h(\psi_h(\bs{s})) ~p(\bs{s})ds$. For example, the density of the QoI
$p_{\Psi_h}(\psi_h)$ can be obtained for
$h(\psi_h(\bs{s})) = \delta(\psi_h-\psi_h(\bs{s}))$
\begin{equation}
p_{\psi_h}(\psi_h)  
=\int \delta(\psi_h -\psi_h(\bs{s}))~ p(\bs{s})\:\dd\bs{s}
\label{eq:ppsih}
\end{equation}
We note at this point that we use upper-case letters to denote the random
variables (e.g.\ $\bs{S}$, $\Psi_h$) and lower-case for the values these take
(e.g.\ $\bs{s}$, $\psi_h$).  Furthermore we use $\psi_h$ and $\psi_h(\bs{s})$
to denote both the values of the QoI as well as the function that provides the
QoI with respect to the input.

The method of choice in high-dimensional settings ($\dim(\bs{S})\gg 1$) is
direct Monte Carlo where the expectation of Equation \eqref{eq:expect} is
estimated by:
\begin{equation}
\hat{\Psi}_h =\frac{1}{N} \sum_{j=1}^N \psi_h(\bs{s}^{j})
\label{eq:mc}
\end{equation}
where $\bs{s}^{j}$ are independent, identically distributed samples drawn from $p(\bs{s})$. The unknown parameters are assumed to depend on a finite number of random variables, with techniques such that the KL theorem discussed above. Hence, a particular grid is introduced in the probability space and for each sample, the solution of a PDE with deterministic coefficients is computed.  The convergence rate is, remarkably independent of the dimension of 
$\bs{S}$ and the error decays as $\mathcal{O}(N^{-1/2})$ \cite{kalos_monte_2008}. Therefore, in problems where each evaluation of $\psi_h(\bs{s})$ poses a 
significant computational burden, the use of direct Monte Carlo can become 
impractical or even infeasible.  

A faster but dimension-dependent method is the Quasi-Monte Carlo Method. It employs low-discrepancy sequences instead of random, or pseudo-random, ones. In this case, sequences of parameters are not randomly chosen, but they are correlated ensuring an order of convergence of $O(N^{-1} (\log\,N)^k)$, where $k$ is the dimensionality of the parameter space.

An improvement over the above Monte Carlo methods is the Stochastic Collocation method \cite{doi:10.1137/100786356,DXiu_JSHesthaven_2005}. The particularity of this method consists in the fact that the choice of Gau\ss-points in the parameter space allows to solve independent problems in space and ensures exponential convergence. A similar approach for solving the stochastic differential equation directly is polynomial chaos (PC) \cite{doi:10.1137/S1064827501387826,doi:10.1142/S0219530512500145}. With this approach, the random process is expanded over a finite number of orthogonal polynomials, such as Hermite basis functions. This method also shows exponential convergence. The application of both KL and PC requires a finite second moment of the random process, while this hypothesis is not necessary for Stochastic Collocation. Clearly, however, these such higher-order methods are still affected by the curse of dimensionality, so in practice is not possible to treat problems with $k > 15$.

\subsection{Multilevel methods}
In order to avoid higher-order methods and deal with the curse of dimensionality, one of the
earliest attempts is to substitute the "true" model with a reduced or surrogate low-fidelity
model~\cite{gian, corrado2015identification}, which we denote by $\psi_l(\bs{s})$. Such an approximation is build
upon the observation that in many applications is clearly possible to
distinguish between \emph{offline} and \emph{online} phases of the workflow,
where the former is typically very expensive and encompasses everything that
can be precomputed in advance, e.g., training a surrogate model before any
patient-specific data becomes available, while the latter is very cheap and
consists only of model evaluations. Unfortunately, this approach has two
significant drawbacks. On the one hand, complex error estimates must be
provided to ensure that the approximation error is within acceptable bounds. On
the other hand, the offline training must be performed only once, which might be impossible in scenarios where the data is changing, such as in the context of precision medicine.

A more straightforward approach to generate low-fidelity models is by coarsening the high-fidelity one. While a single step of coarsening might not reach a desirable reduction of computational cost, the procedure can be iterated several times, obtaining a hierarchy of low-fidelity models. By using a clever partitioning of the samples, most of the simulations can then be performed with the low-resolution models, yielding the so-called multilevel Monte
Carlo method~\cite{gilesmultilevel}. In this case, the direct use of $\psi_h(\bs{s})$ guarantees convergence, while a significant portion of
the computational load is offset to the low-resolution hierarchy. In the case of a 2-level method, the multilevel estimator reads
\begin{equation}
\hat{\Psi}_h =\frac{1}{N_l} \sum_{j=1}^{N_l} \psi_l(\bs{s}^{j}) + \frac{1}{N_h} \sum_{j=1}^{N_h} \left( \psi_h(\bs{s}^{j}) - \psi_l(\bs{s}^{j}) \right)
\label{eq:mlmc}
\end{equation}
Clearly, this approach is convenient if $N_h \ll N_l$ and $\Var[\psi_h(\bs{s}^{j}) - \psi_l(\bs{s}^{j})] \ll \Var[\psi_h(\bs{s}^{j}) ]$. This method extends the idea of control
variables~\cite{fishman_monte_2003} and has shown to be particularly effective in the case of elliptic stochastic partial differential equations~\cite{Harbrecht2013}. However, the creation of the low-resolution model $\psi_l(\bs{s})$ is strictly tied to the possibility to coarsen the geometry of interest, which might not always possible to achieve for real-world problems. This implies that the total number of levels is small, limiting the full potential of this approach~\cite{biehler2015towards}.

\subsection{Multifidelity methods}
In order to produce cheaper low-fidelity models, a more recent idea is to allow for inaccurate models, in the sense of giving up the requirement that their
approximation properties lie within certain error bounds. Instead, they require only
correlation~\cite{peherstorfer2015optimal} or even statistical
dependence~\cite{koutsourelakis2009accurate} with the high-fidelity one. In the former case, the following 2-level estimator can be used
\begin{equation}
\hat{\Psi}_h =\frac{1}{N_h} \sum_{j=1}^{N_h} \psi_h(\bs{s}^{j}) + \alpha \left( \frac{1}{N_l} \sum_{j=1}^{N_l} \psi_l(\bs{s}^{j}) - \frac{1}{N_h} \sum_{j=1}^{N_h} \psi_l(\bs{s}^{j}) \right),
\label{eq:mfmc}
\end{equation}
where the optimal $\alpha$ can be computed from the Pearson correlation coefficient between the two models. This approach has a twofold advantage. On the one hand, it is the statistical
dependence, rather than the error bounds of coarse models, that is crucial to
ensure that propagating uncertainties via the low-fidelity models provides
useful information on the statistics of the high-fidelity quantity-of-interest.
On the other hand, low-fidelity models are not restricted to low-resolution
geometries, therefore solving them can be several orders of magnitude faster
than a coarse model. If the statistical dependence between the models is
present, in practice only $\sim$100 runs of the high-fidelity model are
necessary even for a 2-level method. This number significantly reduces further if a full hierarchy is
produced and exploited~\cite{peherstorfer2015optimal}. Multifidelity methods
have become very popular over the last years and their applications span the
fields of UQ, inverse problems, and optimization~\cite{peherstorfer2016survey}. The key aspect for applying multifidelity methods is therefore the availability
of fast low-fidelity models. Typical examples are the following:
\begin{itemize}
\item
  A mathematical model of expert opinions~\cite{reinert2006including} or
  empirical evidence~\cite{Biehler2017impact}.  In this case, the challenge is
  to create a model that in analytic form can make use of the full stochastic
  input and to quantify its effect on the output.
\item
  Surrogate models, fitting data generated by the high-fidelity model at a
  small number of given input samples \cite{peherstorfer2016survey,
  conti_bayesian_2010}. Unfortunately, these do not perform well in
  high-dimensions of the parameter space ($k \leq 20$), so they are typically augmented by an \emph{a priori}
  dimensionality reduction of the random input vector.  Then, the model is
  trained from this reduced set of variables.
\item
  Projection-based models, such as those based on principal orthogonal
  decomposition~\cite{gian}. This is typically the most involved option, as it
  requires a lot of prior work in deciding which of the techniques available
  are suitable for the problem at hand. In general, all pertinent models and
  training techniques are based on minimizing the difference between the
  full-order and reduced-order model, which is a very strict requirement. In
  the multifidelity approach, it suffices that the outputs are statistically
  dependent.
\item
  Models based on a simplified physical description~\cite{leifsson2010multi}.
  These are based on devising new physical models from first principles, which
  only focus on a portion of the full physical system, or on coarsening the
  existing forward model. This approach has no computational training cost but it could be
  mathematically difficult. The more information  one can introduce, e.g., in
  defining some equivalent/effective properties at the coarsened mesh, the
  better in general the results are.
\end{itemize}

\subsection{Bayesian methods}
While the framework of multifidelity Monte Carlo has proven to greatly reduce
the computing time, it provides point estimates rather than probability
distributions. In fact, a common deficiency of all Monte Carlo-based techniques
is the unavailability of error estimates except for the asymptotic case. Furthermore, they cannot exploit
nonlinear dependencies (which are not reflected in the correlation) between
low- and high-fidelity models in order to accelerate convergence.

\begin{figure}[h]
	\centerline{\includegraphics[width=0.7\textwidth]{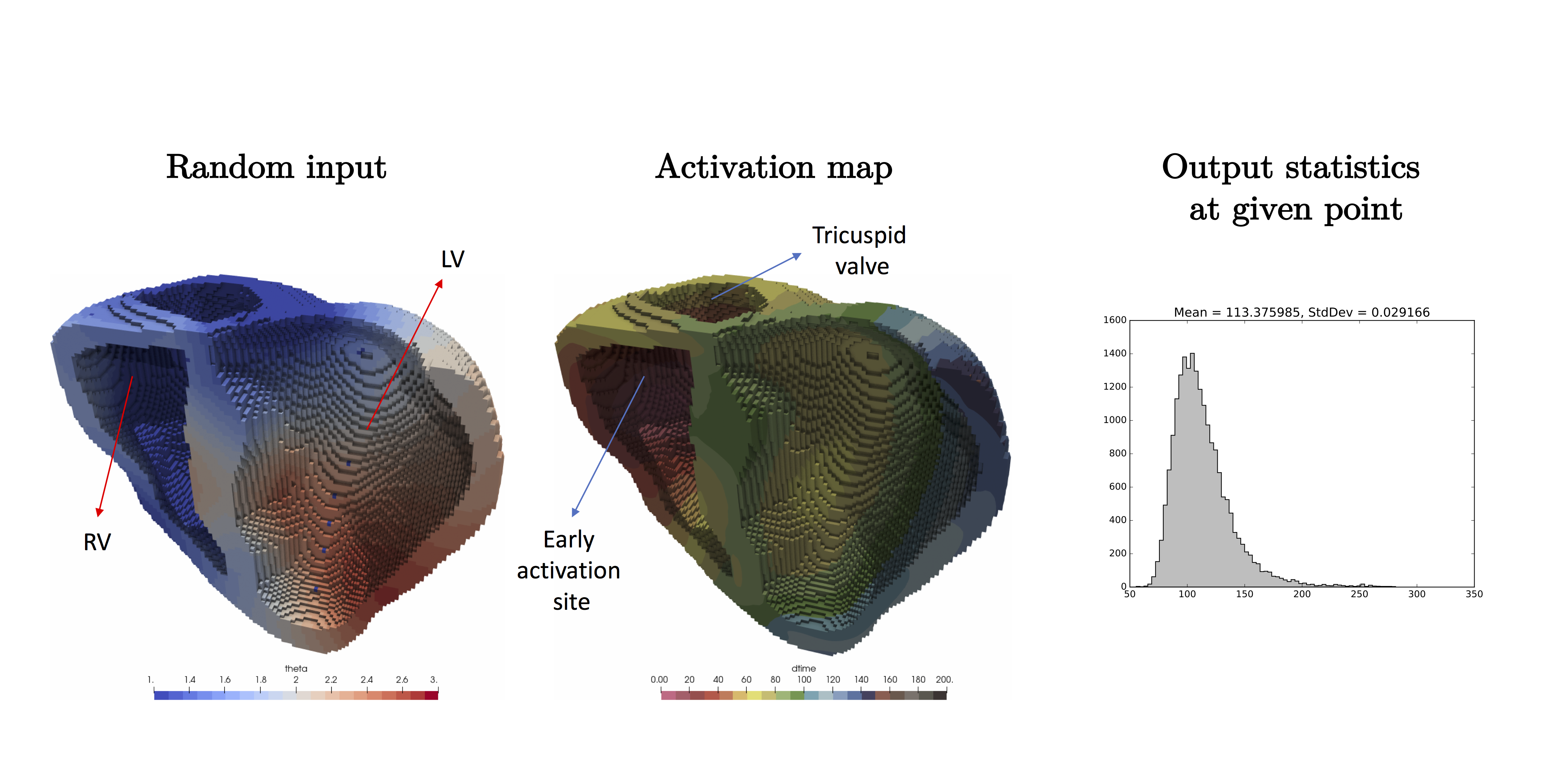}}
	\caption{Uncertainty quantification in cardiac electrophysiology \cite{quaglino2018fast}. Left: one sample of the random perturbation to the fiber field. Center: computed activation map for the random input. Right: probability distribution of the QoI. \label{fig:heart-uq}}
\end{figure}

These issues have motivated several authors to adopt a Bayesian viewpoint~\cite{koutsourelakis2009accurate}, which
automatically augments the estimate with credible intervals and can easily
provide the full probability distribution. This can be used by the analyst to decide whether it is worthwhile to
expend additional effort (in the form of high-order runs)  and can also guide
adaptive enrichment of the number of samples in regions of the random parameter
space that would be most informative for the quantity of interest (QoI). Moreover, it can exploit nonlinear
dependencies between the models, rather than just correlation. This is achieved by fitting a Bayesian regression
between the high-fidelity output and the low-fidelity one. This approach has been successfully employed in \cite{Biehler2017impact,quaglino2018fast} for patient-specific UQ in cardiology, as shown in Figure \ref{fig:heart-uq}.

A Bayesian analog of the frequentist Monte Carlo techniques has been proposed
for general quadrature problems \cite{diaconis_bayesian_1988,
ohagan_monte_1987, ohagan_bayes-hermite_1991, rasmussen_bayesian_2003} and
extended in context of multi-fidelity simulations in
\cite{koutsourelakis2009accurate}.  The fundamental difference of such methods is that the
unknown value of the integral of interest (e.g.\ Equation \eqref{eq:expect})
is treated as a  random variable. The uncertainty in the value
of the integral arises from the fact that it is to expensive to evaluate  the
integrand $ \psi_h(\bs{s})$ at every possible value $\bs{s}$. Instead, the model output is
{\em inferred} by employing a less-expensive,
lower-fidelity model which, in general will not be able to provide the complete picture. This lack
of information translates itself to epistemic uncertainty for the value of the integral
which we attempt to quantify, and potentially reduce by employing more training
data.

Let $\psi_l(\bs{s})$ be the output of a deterministic, lower fidelity and less
expensive model. Given a set  of training data $\mathcal{D}_{N_h} = \{
\psi_l(\bs{s}^{(i)}),\psi_h(\bs{s}^{(i)}) \}_{i=1}^{N_h}$ in the form of low- and
high-fidelity runs, a probabilistic model is built, which is capable of
producing predictive estimates of $\psi_h(\bs{s})$ for any value of $\bs{s}$,
without having to run the expensive high-fidelity model. The premise is that as long as $\psi_l$ exhibits some sort of
statistical dependence with the high-fidelity output $\psi_h$, this could be
accurately learned with a small number of training samples (i.e. high-fidelity
runs) $N_h$. Independently of the dimension of $\bs{s}$, this Bayesian multi-fidelity
strategy advocated implies a $M$-dimensional regression problem, where $M$ is the dimension of the QoI
\cite{koutsourelakis2009accurate}. Therefore, this approach is limited to outputs with a moderate dimensionality. 

Any nonparametric regression tool that 
produces in closed-form expressions for the predictive density $ p(\psi_h | 
\psi_l(\bs{s}), \mathcal{D}_{N_h} )$ can be used. The latter, given the training data $\mathcal{D}_{N_h}$ 
and for each new value $\bs{s}$ of the random input, employs the output of the 
lower-fidelity model  $\psi_l(\bs{s})$ in order to probabilistically infer the value 
$\psi_h$ of the high-fidelity model (for the same $\bs{s}$). The spread of this 
density reflects the epistemic uncertainty in this prediction.
If we denote by $\Psi_{h,\mathcal{D}_{N_h}}(\bs{s})$ the random process that is implied 
by such a probabilistic model and use it in place of the reference $\psi_h(\bs{s})$  
Equation  \eqref{eq:expect}, we  obtain: 
\begin{equation}
  \hat{\Psi}_{h,\mathcal{D}_{N_h}}= \int \Psi_{h,\mathcal{D}_{N_h}}(\bs{s}) ~p(\bs{s})\:\dd\bs{s}
\label{eq:expectbmf}
\end{equation}
 The result of the integral $\hat{\Psi}_{h,\mathcal{D}_{N_h}}$ (even if the 
integration with respect to $\bs{s}$ was performed exactly) would be a random 
variable. Samples of $ \hat{\Psi}_{h,\mathcal{D}_{N_h}}$ are drawn by drawing 
samples of $\Psi_{h,\mathcal{D}_{N_h}}(\bs{s})$.
Furthermore, statistics of $\hat{\Psi}_{h,\mathcal{D}_{N_h}}$ can be computed
\begin{equation}
\begin{split}
\E[\hat{\Psi}_{h,\mathcal{D}_{N_h}} ]
&= \int \E[\Psi_{h,\mathcal{D}_{N_h}}(\bs{s})] p(\bs{s})\dd\bs{s} \\
&= \int \left( \int \psi_h p(\psi_h | \psi_l(\bs{s}), \mathcal{D}_{N_h} ) \dd\psi_h 
\right) p(\bs{s})\:\dd\bs{s}.
\end{split}
\label{eq:meanmean}
\end{equation}
To gain further insight one could compare this with the exact value given in
Equation \ref{eq:expect}. In particular, if we denote  with $p(\psi_h ,
\psi_l)$  the joint density of low- and high-fidelity outputs and given that
the marginal $p(\psi_l)=\int p(\psi_l | \bs{s}) p(\bs{s})\:\dd\bs{s}=\int
\delta(\psi_l -\psi_l(\bs{s})) p(\bs{s})\:\dd\bs{s}$, following standard rules of probability, the following expression is obtained \cite{quaglino2018fast}:
\begin{equation}
\begin{split}
\E[\Psi_h] = \int \biggl(\int \psi_h p(\psi_h | \psi_l(\bs{s}))\dd\psi_h\biggr) p(\bs{s})\:\dd\bs{s}
\end{split}
\label{eq:reexp1}
\end{equation}
By comparing the latter equation with the second line of Equation
\ref{eq:meanmean}, we observe that the exact conditional $ p(\psi_h |
\psi_l(\bs{s}))$ has been replaced by $p(\psi_h | \psi_l(\bs{s}),
\mathcal{D}_{N_h})$ i.e.\ the predictive density implied by the regression model trained on
$\mathcal{D}_{N_h}$. 

 \section{Conclusions}
Multiscale models allow for the treatment of complex phenomena involving different scales, such as, e.g., remodeling and growth of tissues (wound or bone healing, tumor growth, in-stent restenosis), muscular activation, and cardiac electrophysiology. Numerous numerical approaches have been developed to simulate multiscale problems. However, compared to the well-established methods for classical problems, such as solid and fluid mechanics, many questions have yet to be answered in this context. Particularly interesting numerical analysis aspects are 
\begin{itemize}
\itemsep-3pt
\item Consistent definition of a coupled multiscale problem,
\item Define and estimate the error in physical modeling,
\item Properties and stability of the coupled system,
\item Probabilistic treatment of multiscale effects,
\item Parameter estimation,
\item Data-driven simulations.
\end{itemize}

In this paper, we presented an overview of existing models and methods, with particular emphasis on mechanical and bio-mechanical applications. Moreover, we discussed state-of-the-art techniques for multilevel and multifidelity UQ. Finally, we argued that a higher level of abstraction could benefit all fields of multiscale problems, in particular for improving the current approaches of scale generation and fusion. New promising directions are 
\begin{itemize}
\itemsep-3pt
\item Generalize and unify approaches for generating coarse problems,
\item Generalize and unify approaches for scale fusion and transfer,
\item Use modern machine learning and deep learning techniques.
\end{itemize}

In terms of high performance computing, the combination of different simulation tools gives rise to severe technical difficulties, such as the simple exchange of data between different simulation codes, the control of iterative processes, and the introduction of appropriate stopping criteria. Coupled methods have to be developed in order to achieve a good scaling behavior on the upcoming exascale machines. It could be speculated that, in view of the growing core counts, also asynchronous and/or stochastic algorithms might be better suited for the numerical simulation of heterogeneous multiscale problems on the upcoming large scale machines. Such methods are currently not available, at least not to the knowledge of the authors. 

Finally, it can be argued that the huge potential of multiscale simulations can be exploited by intertwining and overcoming the current segregation of modeling, discretizing, solving, and computing solution statistics. The generalization of techniques such as adaptive mesh refinement to {\em model} refinement, where the model is enriched locally, or multigrid methods to incorporate {\em inaccurate} models, as done in multifidelity methods, are examples of desirable yet unexplored avenues of research in this direction.


\bibliographystyle{alpha} 
\bibliography{Krause}

\end{document}